\documentclass[
    ,final            % use final for the camera ready runs
  ]
  {aipproc}

\layoutstyle{6x9}

\usepackage{amsmath}

\def\join{\mathop\vee}
\def\meet{\mathop\wedge}
\def\bor{\mathop{\mathord{\lor}\!\!\!\raise4pt\hbox{$\scriptscriptstyle 2$}\,}}
\def\band{\mathop{\mathord{\land}\!\!\!\lower2pt\hbox{$\scriptscriptstyle 2$}\,}}

\begin{document}

\title{Measuring on Lattices}

\classification{02.10.Ab, 02.10.De, 02.50.Cw}
\keywords      {poset, lattice, algebra, valuation, measure, probability, number theory}

\author{Kevin H. Knuth}{
  address={University at Albany (SUNY), Albany NY, USA}
}

\begin{abstract}
Previous derivations of the sum and product rules of probability theory relied on the algebraic properties of Boolean logic. Here they are derived within a more general framework based on lattice theory. The result is a new foundation of probability theory that encompasses and generalizes both the Cox and Kolmogorov formulations.  In this picture probability is a bi-valuation defined on a lattice of statements that quantifies the degree to which one statement implies another.  The \emph{sum rule} is a constraint equation that ensures that valuations are assigned so as to not violate associativity of the lattice join and meet.  The product rule is much more interesting in that there are actually \emph{two product rules}: one is a constraint equation arises from associativity of the direct products of lattices, and the other a constraint equation derived from associativity of changes of context.  The generality of this formalism enables one to derive the traditionally assumed condition of additivity in measure theory, as well introduce a general notion of product.  To illustrate the generic utility of this novel lattice-theoretic foundation of measure, the sum and product rules are applied to number theory.  Further application of these concepts to understand the foundation of quantum mechanics is described in a joint paper in this proceedings.
\end{abstract}

\maketitle

%%%%%%%%%%%%%%%%%%%%%%%%%%%%%%%%%%%%%%%%%%%%
%% MAINMATTER
%%%%%%%%%%%%%%%%%%%%%%%%%%%%%%%%%%%%%%%%%%%%

\section{Introduction}

A lattice is an algebra.  Where an algebra considers a set of elements along with a set of operations that takes one or more elements to another element, the lattice considers a set of elements along with a binary ordering relation that sets up a hierarchy among the elements.  The algebraic perspective is operational, whereas the lattice perspective is structural.  Both the operational and structural relationships among elements are useful.

\subsection{Posets, Lattices and Algebras}
\begin{figure}[t]
  \label{fig:posets}
  \includegraphics[height=.15\textheight]{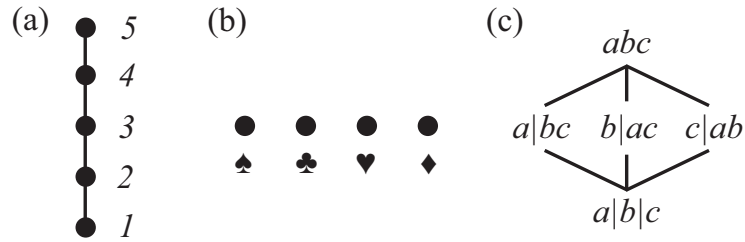}
  \caption{Three basic examples of posets. (a) The integers ordered by the usual $\leq$ form a \emph{chain}.  The element $2$ is drawn above $1$ since $1 \leq 2$, and they are connected by a line because $2$ covers $1$ in the sense that there is no integer $x$ between $2$ and $1$ such that $1 \leq x \leq 2$.  (b) The four card suits are incomparable under a wide variety of card game rules and we draw them side-by-side to express this.  This configuration is called an \emph{antichain}.  (c) The set of partitions of three elements $a$, $b$ and $c$ ordered by partition containment forms a more complex poset that exhibits both chain and antichain behavior.  One chain consists of the elements $a | b | c$, $a | bc$, and $abc$ since each successive partition contains the previous.  The elements $a | bc$, $b | ac$, and $c | ab$ form an antichain because not one of these three partitions contains another.}
\end{figure}

Two elements of a set are ordered by comparing them according to a binary ordering relation, generically denoted $\leq$ and read `\emph{is included by}'.  One of the first examples that comes to mind may be the ordering of the integers according to the usual meaning of the symbol $\leq$ `\emph{is less than or equal to}'.  The ordering that results is called a \emph{chain} (Fig. \ref{fig:posets}a).  To illustrate the hierarchy, we simply draw element $B$ above element $A$ if $A \leq B$ and connect them with a line if there does not exist an element $X$ in the set such that $A \leq X \leq B$.

In some cases, elements of the set are incomparable to one another, as in the common example of comparing apples and oranges.  Figure \ref{fig:posets}b shows an \emph{antichain} of card suits where the elements are placed side-by-side to indicate that no element includes any other.

More interesting examples involve both inclusion and incomparability, which is why these structures generally are called \emph{partially ordered sets}, or \emph{posets} for short.  Figure \ref{fig:posets}c illustrates the lattice that results from partitioning three objects.  One could consider all three objects together $a b c$, or each separately $a | b | c$.  We could also partition the objects in these ways: $a | b c$, $b | a c$ or $c | a b$.  These partitions can be compared according to a relation that decides whether one partition includes another.  The partition $a b c$ includes the partition $a | b | c$ since it can be obtained by simply sub-dividing $abc$.  However, the partitions $c | a b$ and $a | b c$ are incomparable since, for example, there is no way to sub-divide the partition $c | a b$ to obtain $a | b c$.

\begin{figure}[t]
  \label{fig:lattice}
  \includegraphics[height=.12\textheight]{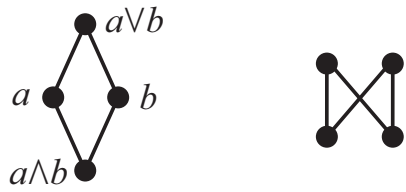}
  \caption{The poset on the left is a simple lattice, which illustrates the join $\join$ and the meet $\meet$.  The poset on the right is not a lattice since the pair of elements on the bottom do not have a unique least upper bound. Similarly, the pair of elements at the top do not have a unique greatest lower bound.}
\end{figure}

Given a set of elements in a poset, their \emph{upper bound} is the set of elements that contain them.  For example, the upper bound of the partition $c | a b$ in Fig. \ref{fig:posets}c is the set $\{a b c\}$.  Given a pair of elements $x$ and $y$, the least element of their upper bound is called the \emph{join}, denoted $x \vee y$.  The \emph{lower bound} of a pair of elements is defined dually by considering all the elements that the pair of elements contain.  The greatest element of the lower bound is called the \emph{meet}, denoted $x \wedge y$. A \emph{lattice} is a partially ordered set where each pair of elements has a unique meet and a unique join (Fig. \ref{fig:lattice}).  Graphically, the join can be found by starting at both elements and following the lines upward until they first intersect.  The meet is found similarly by moving downward.  There often exist elements that are not formed from the join of any pair of elements.  These elements are called \emph{join-irreducible elements}.  \emph{Meet-irreducible elements} are defined similarly. For example, the partitions $a | b c$, $b | a c$ or $c | a b$ cannot be formed by joining any other pair of partitions and therefore are join-irreducible.  In this case, these elements are also meet-irreducible.

We can choose to view the join and meet as algebraic operations that take any two lattice elements to a unique third lattice element.  From this perspective, the lattice is an algebra.  This gives us both a structural and operational perspective which are related by a set of equations called \emph{consistency relations}
\begin{equation}
\label{eq:lattice-algebra} x \leq y \qquad\Longleftrightarrow\qquad
   \begin{array}{rl}
       x \vee y = y \\
       x \wedge y = x \end{array}
\end{equation}

Given a specific lattice, we find that the consistency relations result in a specific algebraic identity.
For example, the integers ordered by the usual `\emph{less than or equal to}' leads to
\begin{equation}
\label{eq:leq} x \leq y \qquad\Longleftrightarrow\qquad
   \begin{array}{rl}
       \max(x, y) = y \\
       \min(x, y) = x \end{array}
\end{equation}
whereas the positive integers ordered by `\emph{divides}' leads to
\begin{equation}
\label{eq:divides} x \mid y \qquad\Longleftrightarrow\qquad
   \begin{array}{rl}
       \mathrm{lcm}(x, y) = y \\
       \gcd(x, y) = x \end{array}
\end{equation}
Sets ordered by the usual `\emph{is a subset of}' leads to
\begin{equation}
\label{eq:subseteq} x \subseteq y \qquad\Longleftrightarrow\qquad
   \begin{array}{rl}
       x \cup y = y \\
       x \cap y = x \end{array}
\end{equation}

\section{Spaces}
Lattices enable us to describe our world and encode the order that we observe.  However, a description of the world is not the same as a description of our state of knowledge about the world.  Here we examine the different spaces that we can construct using lattices.

\subsection{Components and State Space}
\begin{figure}
  \label{fig:states}
  \includegraphics[height=.15\textheight]{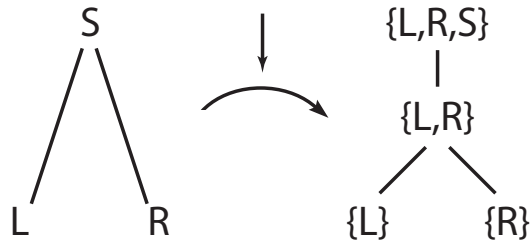}
  \caption{The poset on the left encodes the relationship among the components of a bridge: the left side $L$, the right side $R$, and the span $S$, which sits on top.  By taking the set of downsets ordered by set inclusion, we obtain a lattice describing all the possible states of a bridge where each state is a set of bridge components.}
\end{figure}

Systems are often made of components, and these components may have a natural order.  Consider the components of a bridge: the left side $L$, the right side $R$, and the span $S$, which sits on top of the two sides.  The relationship between the components of a bridge is illustrated in Fig. \ref{fig:states}.

The states of a bridge can be constructed from the components by taking the set of \emph{downsets} of bridge components and ordering them according to set inclusion.  A downset is a set of components, such that if an element is a member of the set, then all elements that it includes are also members of the set.  The result is a lattice of possible bridge states where each state is a set of bridge components.  Note that the downset construction prevents one from including sets of components such as $\{L, S\}$  from being included as a state, since the span can only be present if both the left and right sides are present.

Not all state spaces are constructed from components.  Instead the states may represent the most elementary description of the system. Often, the poset of states forms an antichain consisting of independent and mutually exclusive elements.

\begin{figure}[b]
  \label{fig:statements}
  \includegraphics[height=.175\textheight]{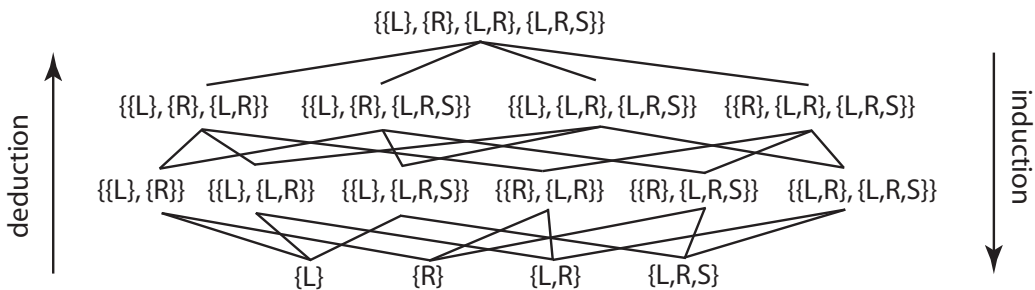}
  \caption{This figure depicts the lattice of all statements that can be made about the state of the bridge.  A statement is defined as a set of potential system states.  The ordering relation of set inclusion naturally encodes logical implication (deduction), such that a statement implies all the statements above it.}
\end{figure}

\subsection{Hypothesis Space}
The bridge may be in one of four possible states of construction.  A given individual may not know precisely which state the bridge is in, but may have some information that rules out some states, but not others.  This set of \emph{potential states} defines what one can say about the state of the bridge. For this reason, I call a set of potential states a \emph{statement}.  In this sense, a statement describes a state of knowledge about the state of the bridge. Note that quantification will come later.

The lattice of statements is generated by taking the powerset, which is the set of all possible subsets of the set of all states, and ordering them according to set inclusion.  In this example, since states are sets of components, statements are sets of sets.  Other examples have been given in previous works \cite{Knuth:laws,Knuth:me08}. Given that there are four possible states, there are $2^4 = 16$ statements including the null set, which is generically called the \emph{bottom} and denoted $\bot$ since it resides at the bottom of the lattice. The bottom element is often omitted from the diagram due to the fact that it represents the logical absurdity.  These statements, illustrated in Fig. \ref{fig:statements}, form the \emph{hypothesis space}.  The ordering relation of set inclusion naturally encodes logical implication, such that a statement implies all the statements above it.

The statements along the bottom of the lattice represent states of knowledge where the state of the bridge is known with certainty.  The statement at the top is the \emph{truism}, generically called the \emph{top} and denoted $\top$, which represents the state of knowledge where one only knows that the bridge can be in one of four possible states.  Intermediate statements represent intermediate states of knowledge.  For example, you may know that the right side $R$ will be constructed by either Worker A, who is a very hard worker, or Worker B who is very lazy.  In this case your state of knowledge about the state of the bridge would probably be the intermediate state $\{\{L\}, \{L,R,S\}\}$, which says that either just the left side is up or the whole bridge is finished.

\emph{Logical deduction} is straightforward in this framework since a statement in the lattice implies (is included by) every statement above it with certainty.  \emph{Logical induction} works backwards.  One would like to quantify the \emph{degree} to which one's current state of knowledge implies a statement of greater certainty \emph{below} it.  Since statements do not imply statements below them, this requires a generalization of the algebra representing the ordering.  In the next section we will generalize this algebra to a calculus by introducing quantification.  The result is a measure, called \emph{probability}, that quantifies the degree to which one statement implies another.

\subsection{Inquiry Space}
One can take this idea of states and statements further.  By defining a question in terms of the statements that answer it \cite{Cox:1979}, one can generate the lattice of questions by taking downsets of statements \cite{Knuth:Questions}. That is, if a given statement answers a question, then all the statements that imply that statement also answer the question.  In the case of our bridge, the lattice of questions has 167 elements.

Just as some statements imply other statements, some questions answer other questions.  Answering a given question in the lattice will guarantee that you have answered all the questions above it.  One can also generalize this algebra to a calculus by introducing quantification.  The result is the \emph{inquiry calculus}, which is based on a measure, called \emph{relevance}, that quantifies the degree to which one question answers another.

For more details on the lattice of questions and inquiry, see the following references \cite{Knuth:Questions, Knuth:duality,Knuth:WCCI06,Knuth:me08}, as well as the paper by Julian Center in this proceedings \cite{Center:Inquiry}.

\section{Quantification}
An algebra can be extended to a calculus by defining functions that take lattice elements to real numbers.  This enables one to \emph{quantify} the relationships between the lattice elements.  A valuation $v$ is a function that takes a single lattice element $x \in L$ to a real number $v(x)$ in a way that respects the partial order, so that $v(x) \leq v(y)$ iff $x \leq y$.  This means that the lattice structure imposes constraints on the valuation assignments, which can be expressed as a set of constraint equations.

The valuation assigned to element $x$ can be defined with respect to a second lattice element $y$ called the \emph{context}.  The result is a function called a bi-valuation $w(x \mid y) = v_y(x)$, which takes two lattice elements $x$ and $y$ to a real number.  Here a solidus is used as an argument separator so that one reads $w(x \mid y)$ as the degree to which $y$ includes $x$.

In the following sections, we consider three symmetries and identify the constraints that they impose on valuation and bi-valuation assignments.  The first two symmetries result from the lattice structure and thus impose the same constraints on the valuation and bi-valuation assignments; whereas the last symmetry considered is specific to bi-valuations.

\begin{figure}
  \label{fig:sum}
  \includegraphics[height=.15\textheight]{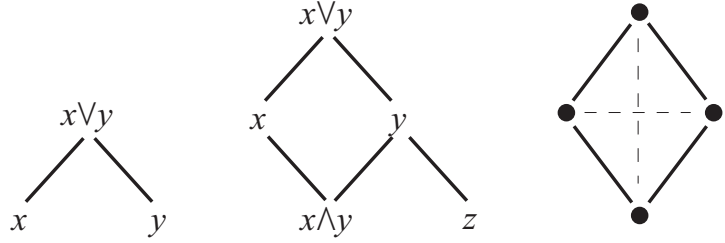}
  \caption{The poset on the left is used to establish the additive nature of the valuation.  The poset in the center is used to establish the sum rule for the lattice in general.  The cartoon on the right illustrates the symmetry of the sum rule.  The sum of the valuations of the elements at the top and bottom of the diamond equals the sum of the valuations of the elements on the right and left sides.  These dashed lines conveniently form a plus sign reminding us of the sum rule.}
\end{figure}

\subsection{Sum Rule}
We begin by considering a special case depicted in Fig. \ref{fig:sum} (left) of two elements $x$ and $y$ with join $x \join y$ and a null meet $x \meet y = \bot$ (not shown).  The value we assign to the join $x \join y$, written $u(x \join y)$, must be a function of the values we assign to both $x$ and $y$, $u(x)$ and $u(y)$, since if there did not exist any functional relationship, then the valuation could not possibly reflect the underlying lattice structure.  We write this functional relationship in terms of an unknown binary operator $\oplus$
\begin{equation}
u(x \join y) = u(x) \oplus u(y).
\end{equation}
We now consider another case where we have three elements $x$, $y$, and $z$, such that their meets are again disjoint.  The least upper bound of these three elements can be written in these two different ways: $x \join (y \join z)$ and $(x \join y) \join z$.  The value we assign to this join can also be written in two different ways
\begin{equation}
u(x) \oplus \big(u(y) \oplus u(z)\big) = \big(u(x) \oplus u(y)\big) \oplus u(z).
\end{equation}
This is a functional equation for the operator $\oplus$ for which the general solution is given by Aczel \cite{Aczel:FunctEqns}
\begin{equation}
f(u(x \join y)) = f(u(x)) + f(u(y)),
\end{equation}
where $f$ is an arbitrary invertible function, so that many valuations are possible.
We take advantage of this freedom to choose a valuation $v(x) = f(u(x))$ that simplifies this constraint
\begin{equation} \label{eq:simple-sum}
v(x \join y) = v(x) + v(y).
\end{equation}
By letting $x = \bot$, equation (\ref{eq:simple-sum}) implies that $v(\bot) = 0$.

Now that we have a constraint on the valuation for our simple example, we seek the general solution for the entire lattice.  To derive the general case, we consider the lattice in Figure \ref{fig:sum} (center) and note that the elements $x \meet y$ and $z$ have a null meet, as do the elements $x$ and $z$.  Applying (\ref{eq:simple-sum}) to these two cases, we get
\begin{eqnarray}
v(y) & = & v(x \meet y) + v(z)\\
v(x \join y) & = & v(x) + v(z)
\end{eqnarray}
Simple substitution results in the general constraint equation known as the \emph{sum rule}
\begin{equation}
v(x \join y) = v(x) + v(y) - v(x \meet y).
\end{equation}
In general for bi-valuations we have
\begin{equation}
v(x \join y \mid t) = v(x \mid t) + v(y \mid t) - v(x \meet y \mid t).
\end{equation}
for any context $t$. Note that the sum rule is not focused solely on joins since it is symmetric with respect to interchange of joins and meets.

At this point, we have \emph{derived} additivity of the measure, which is considered to be an \emph{axiom} of measure theory.  This is significant in that associativity constrains us to have additive measures---there is no other option.  The cartoon at the right of Fig. \ref{fig:sum} illustrates the symmetry of the sum rule.  The sum of the valuations of the elements at the top and bottom of the diamond equals the sum of the valuations of the elements on the right and left sides
\begin{equation} \label{eq:sum-rule}
v(x \join y) + v(x \meet y) = v(x) + v(y).
\end{equation}

\subsection{Lattice Products}
Given the linearity of the constraint imposed by associativity (\ref{eq:sum-rule}), the only remaining freedom is that of rescaling.  This means that any further constraints must have a multiplicative form.  One can combine two lattices via the lattice product where elements are combined in as in a Cartesian product.  That is, the product of a lattice $X$ with a lattice $Y$ will result in a lattice $X \times Y$ with elements of the form $(x, y)$, where $x \in X$ and $y \in Y$.  The lattice product is associative, so that for three lattices $X$, $Y$, and $Z$, we have
\begin{equation}
(X \times Y) \times Z = X \times (Y \times Z)
\end{equation}
with elements of the form $(x,y,z)$.

The valuation assigned to an element $(x, y)$ clearly must be a function of the valuations assigned to $x$ and $y$.  However, the associative relation above imposes a constraint on the valuations assigned to the elements of the lattice product.  This constraint ultimately results in the following relation
\begin{equation}
v((x, y)) = v(x) v(y),
\end{equation}
which is a \emph{product rule} that applies when one combines two independent spaces.

\begin{figure}
  \label{fig:prod}
  \includegraphics[height=.25\textheight]{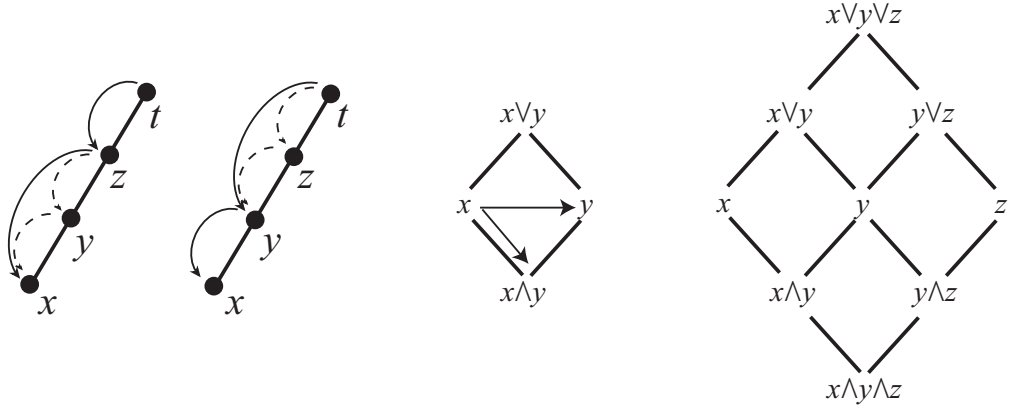}
  \caption{The chains on the left illustrate the changes in context used to derive the chain rule.  The diamond in the center illustrates that the degree to which $x$ includes $x \meet y$ equals the degree to which $x$ includes $y$.  The lattice on the right is used to derive the general product rule describing context change.  See text below for details.}
\end{figure}

\subsection{The chain rule}
We now focus on bi-valuations and explore changes in context. We begin with a special case and consider four ordered elements $x \le y \le z \le t$.  The relationship $x \le z$ can be divided into two relations, $x \le y$ and $y \le z$.  In the event that $z$ is considered to be the context, this sub-division implies that the context can be considered in parts.  Thus the bi-valuation we assign to $x$ with respect to context $z$, $w(x \mid z)$, must be related to both the bi-valuation we assign to $x$ with respect to context $y$, $w(x \mid y)$, and the bi-valuation we assign to $y$ with respect to context $z$, $w(y \mid z)$.  That is, there exists a binary operator $\odot$ that relates the bi-valuations assigned to the two steps to the bi-valuation assigned to the one step
\begin{equation}
    w(x\mid z) = w(x\mid y) \odot w(y\mid z)\,.
\end{equation}
By extending this to three steps (Fig. \ref{fig:prod}, left) and considering the bi-valuation $w(x\mid t)$ relating $x$ and $t$, via intermediate contexts $y$ and $z$, results in another associativity relationship
\begin{equation}
    \big( w(x\mid y) \odot w(y\mid z) \big) \odot w(z \mid t) = w(x\mid y) \odot \big( w(y\mid z) \odot w(z \mid t) \big)
\end{equation}
Using the associativity theorem again results in a constraint equation for non-negative bi-valuations involving changes in context \cite{Skilling:me08}.  We call this the \emph{chain rule}
\begin{equation}
w(x \mid z) = w(x \mid y) w(y \mid z)\,.
\end{equation}

\subsection{The product rule}
We now focus on extending these results to the more general case illustrated with the lattice on the right side of Fig. \ref{fig:prod}.  To begin we focus on the small diamond defined by $x$, $y$, $x \lor y$, and $x \land y$.  If we consider the context to be $x$, the sum rule gives
\begin{equation}
w(x\mid x) + w(y\mid x) = w(x \lor y\mid x) + w(x \land y\mid x).
\end{equation}
Since $x \leq x$ and $x \leq x \lor y$, we have $w(x\mid x) = w(x \lor y\mid x) = 1$, and the sum rule reduces to
\begin{equation}\label{eq:prodstep1}
w(y\mid x) = w(x \land y\mid x).
\end{equation}
This relationship is illustrated by the equivalence of the arrows in the diamond in the center of Fig. \ref{fig:prod}.
This result will used several times in the derivation that follows.

Consider the chain where the bi-valuation $w(x \meet y \meet z \mid x)$ with context $x$ is decomposed into two parts by introducing the intermediate context $x \meet y$.  The chain rule gives
\begin{equation}\label{eq:prodchain}
w(x \meet y \meet z\mid x) = w(x \meet y \meet z \mid x \meet y)\,w(x \meet y \mid x).
\end{equation}
To simplify this relation, consider the diamond defined by $x \meet y \meet z$, $x \meet y$, $y \join z$, $z$ to obtain
\begin{equation}\label{eq:prodstep2}
w(x \meet y \meet z \mid x \meet y) = w(z \mid x \meet y).
\end{equation}
Similarly, consider the diamond defined by $x$, $x \join y$, $y \meet z$, and $x \meet y \meet z$ to obtain
\begin{equation}\label{eq:prodstep3}
w(x \meet y \meet z \mid x) = w(y \meet z \mid x).
\end{equation}
Substituting (\ref{eq:prodstep1}),(\ref{eq:prodstep2}), and (\ref{eq:prodstep3}) into (\ref{eq:prodchain}) results in
\begin{equation}
w(y \meet z \mid x) = w(z \mid x \meet y)\,w(y \mid x),
\end{equation}
which is the general \emph{product rule} for context change.

\section{The valuation calculus}
We have derived that associativity of the lattice join leads to the sum rule for valuations
\begin{equation}
    v(x \join y) + v(x \meet y) = v(x) + v(y)\,.
\end{equation}
which is a key axiom of measure theory.  Associativity of the lattice product imposes an additional constraint, which results in a product rule
\begin{equation}
    v((x, y)) = v(x)v(y)\,.
\end{equation}

Extending the concept of valuation to that of a context-dependent bi-valuation, we obtain a sum rule
\begin{equation}
    w(x \join y \mid t) + w(x \meet y \mid t) = w(x \mid t) + w(y \mid t)\,,
\end{equation}
a product rule for combining spaces
\begin{equation}
    w((x, y) \mid (t_x, t_y)) = w(x \mid t_x)w(y \mid t_y)\,,
\end{equation}
and a product rule for context change
\begin{equation}
w(y \meet z \mid x) = w(z \mid x \meet y)\,w(y \mid x)\,.
\end{equation}
This calculus of valuations not only derives the condition of additivity of measures, but also generalizes traditional measure theory in an important way by introducing the concept of context.  In addition to the sum rule of measure theory, composition of spaces and changes in context introduce two distinct product rules.

\section{Applications}
Since these results are valid for \emph{all} lattices, we can apply them to a wide array of applications.  Applying the valuation calculus to the lattice of statements results in probability theory where the bi-valuations are conditional probabilities.  Bayes theorem can be seen to implement a particular change in context.  These results can also be applied to the lattice of questions resulting in the inquiry calculus, which is an analogous calculus of a measure on questions called relevance.

These results can also be applied to the lattice of experimental setups in quantum mechanics with a significant difference where pairs of real numbers are needed to quantify the lattice elements rather than scalars as in the examples considered here.  Joint papers written with Philip Goyal and John Skilling in this proceedings \cite{GKS:me09} and elsewhere \cite{GKS:PRA} show why quantum amplitudes behave as complex numbers obeying Feynman's rules using methods very similar to those described here.

To reinforce the idea that these sum and product rules are not limited to the domain of probability and are instead applicable to any consistent valuation over all lattices, we will briefly explore application of the valuation calculus to an unlikely domain---number theory.

\begin{figure}
  \label{fig:number-theory}
  \includegraphics[height=.2\textheight]{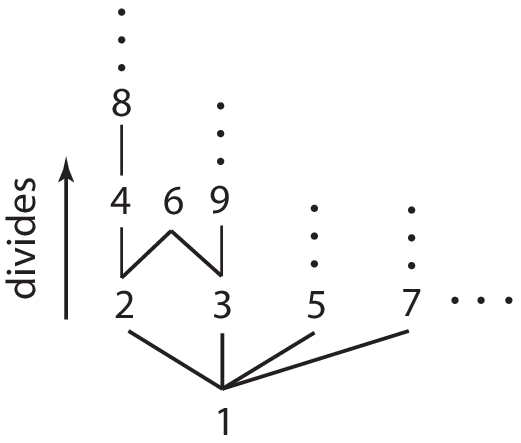}
  \caption{The lattice generated by ordering the integers according to whether one integer divides another.  The bottom element is unity. The primes are the atomic elements, and the primes and powers of primes are the join-irreducible elements.  For the purposes of quantification, it is assumed that the lattice is not infinite in extent so that there exists a greatest integer in the set. The join of two elements is their least common multiple, and the meet is their greatest common divisor.}
\end{figure}

\subsection{Number Theory and Sum and Product Rules}
We consider the lattice of integers less than a greatest integer $t$ ordered by division.  A portion of the lattice is illustrated in Figure \ref{fig:number-theory}.  The bottom element is unity. The primes are the atomic elements, and the primes and powers of primes are the join-irreducible elements, which provides some insight as to why they have a special status in number theory.  The consistency relation (\ref{eq:divides}) illustrates that the join and meet are the number theoretic operations least common multiple ($\mathrm{lcm}$) and greatest common divisor ($\gcd$), respectively.  Since the operation $\mathrm{lcm}$ distributes over $\gcd$ and vice versa, the lattice is a distributive lattice.  This means that we have the freedom to assign valuations to the join-irreducible elements (primes and powers of primes) arbitrarily.

An interesting valuation assignment is one where $v(p) = \log p$ for all $p$ such that $p$ is a prime or power of a prime.  For two primes $p$ and $q$, we have according to (\ref{eq:simple-sum})
\begin{equation}
v(\mathrm{lcm}(p, q)) = \log p + \log q.
\end{equation}
Specifically, for $p = 2$ and $q = 3$,  we have that
\begin{equation}
v(6) = \log 2 + \log 3 = \log 6.
\end{equation}
In general, the sum rule holds.  For example,
\begin{equation}
v(12) = \log 4 + \log 6 - \log 2 = \log 12,
\end{equation}
since $\mathrm{lcm}(4, 6) = 12$ and $\gcd(4, 6) = 2$.

The chain rule and product rule hold as well.  Let the greatest integer in the set be denoted by $t$, and assign $w(p \mid t) = \log p$ for all $p$ such that $p$ is a prime or power of a prime.  It is straightforward to show that $d(n \mid 1) = 1$ so that the degree to which unity divides any integer $n$ is 1.  Similarly, $d(n \mid n) = 1$ and $d(n \mid m) = 1$ for all $m \leq n$.

We can use the product rule to write
\begin{eqnarray}
d(m \meet n \mid t) & = & d(m \mid t) d(n \mid m \meet t) \\
& = & d(n \mid t) d(m \mid n \meet t)
\end{eqnarray}
and equate the two right-hand sides as in the derivation of Bayes' Theorem to obtain
\begin{equation}
d(m \mid n \meet t) = \frac{d(m \mid t) d(n \mid m \meet t)}{d(n \mid t)}.
\end{equation}
In the case where we let $m \leq n < t$ we have
\begin{equation}
d(m \mid n) = \frac{\log m}{\log n},
\end{equation}
which for $m = 2$ and $n = 4$ gives us the degree to which four divides two
\begin{equation}
d(2 \mid 4) = \frac{\log 2}{\log 4} = 1/2.
\end{equation}

\section{Conclusion}
In this paper we derive the valuation calculus. Associativity of the join gives rise to the sum rule, which is symmetric with respect to interchange of joins and meets.  Associativity of the lattice product results a product rule, which dictates how valuations are to be combined when taking lattice products.  Associativity of changes of context result in a product rule for bi-valuations that dictate how valuations should be manipulated when changing context.  These results are valid for all lattices---not just the Boolean algebra of probability theory.  This new theory of quantification opens the door to a wide variety of novel applications, such as decision theory and concept lattices.

With respect to probability theory the result of this work is a new foundation that encompasses and generalizes both the Cox and Kolmogorov formulations.  By introducing probability as a bi-valuation defined on a lattice of statements we can quantify the degree to which one statement implies another.  This generalization from logical implication to degrees of implication not only mirrors Cox's notion of plausibility as a degree of belief, but includes it.  The main difference is that Cox's formulation is based on a set of desiderata derived from his particular notion of plausibility; whereas here the symmetries of lattices in general form the basis of the theory and the \emph{meaning} of the derived measure is inherited from the ordering relation, which in this case is implication.  The fact that these lattices are derived from sets means that this work encompasses Kolmogorov's formulation of probability theory as a measure on sets.  However, mathematically this theory improves on Kolmogorov's foundation by \emph{deriving}, rather than assuming, summation.  Furthermore, this foundation further extends Kolmogorov's measure-theoretic foundation by introducing the concept of context.  This leads directly to probability necessarily being conditional, and Bayes' Theorem follows as a direct result of the product rule by which it implements change in context.

By better understanding the foundation, we expect to be able to more readily extend our theory of inference into greater domains.  This will include a better understanding of maximum entropy and perhaps lead to entropic inference where both data and constraints are readily combined to aid our inferences.  In another direction we consider inference in quantum mechanics, and our first steps in this direction can be found in a joint paper in this volume \cite{GKS:me09} and elsewhere \cite{GKS:PRA}.

%%%%%%%%%%%%%%%%%%%%%%%%%%%%%%%%%%%%%%%%%%%%%%%%
%% BACKMATTER
%%%%%%%%%%%%%%%%%%%%%%%%%%%%%%%%%%%%%%%%%%%%%%%%

\begin{theacknowledgments}
I would like to thank John Skilling for his encouragement and strong support of this work, as well the new ideas and directions that he has introduced.  I would also like to thank Janos Acz\'{e}l, Ariel Caticha, Julian Center, Philip Goyal, Steve Gull, Jeffrey Jewell, Vassilis Kaburlasos, and Carlos Rodr\'{i}guez for inspiring discussions, invaluable remarks and comments, and much encouragement.  A special thanks goes to Tom Loredo for a careful reading of this manuscript and the valuable comments he provided.  This work was supported in part by the College of Arts and Sciences and the College of Computing and Information of the University at Albany (SUNY), the NASA Applied Information Systems Research Program (NASA NNG06GI17G) and the NASA Applied Information Systems Technology Program (NASA NNX07AD97A).
\end{theacknowledgments}



\bibliographystyle{aipproc}   % if natbib is available

\bibliography{knuth}

\IfFileExists{\jobname.bbl}{}
 {\typeout{}
  \typeout{******************************************}
  \typeout{** Please run "bibtex \jobname" to optain}
  \typeout{** the bibliography and then re-run LaTeX}
  \typeout{** twice to fix the references!}
  \typeout{******************************************}
  \typeout{}
 }

\end{document}